\documentclass[11pt]{amsart}

 \usepackage{amsmath,amsthm,amsfonts,amssymb,verbatim}
 \usepackage{url}
 \usepackage{graphicx}
 \usepackage[all]{xy}
 \usepackage{wrapfig}
 \usepackage{picins}
 \usepackage{pinlabel}
 \usepackage{subfigure}
 
\let\proved\qedhere 
 
\newcommand{\Br}{\mathbf{\Sigma}}

\newcommand{\bZ}{\mathbb{Z}}

\newcommand{\bF}{\mathbb{F}}

\newcommand{\rs}{\textstyle\frac{r}{s}}

\newcommand{\overzero}{\textstyle\frac{1}{0}}

\newcommand{\fibre}{\varphi}

\newcommand{\into}{\hookrightarrow}

\newcommand{\Khred}{\widetilde{\operatorname{Kh}}}
\newcommand{\rk}{\operatorname{rk}}

\newtheorem{theorem}{Theorem}

\newtheorem{proposition}[theorem]{Proposition}

\newtheorem*{namedtheorem}{\theoremname}
\newcommand{\theoremname}{testing}
\newenvironment{named}[1]{\renewcommand{\theoremname}{#1}
        \begin{namedtheorem}}
        {\end{namedtheorem}}
        
\title[Khovanov homology and two-fold branched covers]{A remark on Khovanov homology and two-fold branched covers}
\date{July 13, 2009}

\author[Liam Watson]{Liam Watson}
\thanks{Supported by a Canada Graduate Scholarship (NSERC)}
\address{D\'epartement de Math\'ematiques, Universit\'e du Qu\'ebec \`a Montr\'eal,  Montr\'eal Canada.}
\email{liam.watson@cirget.ca}
\urladdr{http://www.cirget.uqam.ca/~liam}

\begin{document}

\begin{abstract}
Examples of knots and links distinguished by the total rank of their Khovanov homology but sharing the same two-fold branched cover are given. As a result, Khovanov homology does not yield an invariant of two-fold branched covers. 
\end{abstract}

\maketitle

Mutation provides an easy method for producing distinct knots sharing a common two-fold branched cover: the mutation in the branch set corresponds to a trivial surgery in the cover. Due to a result of Wehrli \cite{Wehrli2009,Wehrli2007} (see also Bloom \cite{Bloom2009}), this provides a range of examples of manifolds that branch cover $S^3$ in more than one way, but for which the distinct branch sets have identical rank in their respective Khovanov homology groups over $\bF=\bZ/2\bZ$. 

From this point of view this fact is not completely surprising, as Khovanov homology is closely related to the Heegaard-Floer homology of two-fold branched covers \cite{OSz2005-branch}. Indeed, this is made precise in Bloom's proof of mutation invariance \cite{Bloom2009}. More generally however, the following question has been posed by Ozsv\'ath \cite{Lipshitz2008}: is Khovanov homology an invariant of the two-fold branched cover? To be precise, this question asks if the {\em total rank} of the reduced Khovanov homology (over  $\bF$) might be an invariant of two-fold branched covers. This arises naturally in relation to the open problem of extending Khovanov's invariant to manifolds other than $S^3$ \cite{Lipshitz2008,Ozsvath2008}.  

This short note gives a negative answer.  

\begin{named}{Theorem}
The total rank of Khovanov homology is not an invariant of two-fold branched covers. \end{named}

This theorem is proved by example: we construct manifolds that are two-fold branched covers of $S^3$ in two different ways, and for which the pair of branch sets is distinguished by the total rank in Khovanov homology. We work with the reduced version of Khovanov homology, denoted $\Khred$, with $\bF$ coefficients \cite{Khovanov2000,Khovanov2003}. 

\subsection*{Surgery on torus knots}

Let $S^3_{r/s}(K)$ denote the result of $\rs$-framed surgery on a knot $K\into S^3$, and let $T_{p,q}$ denote the positive $(p,q)$ torus knot in $S^3$ (with $0<p<q$). Note that, as we will only consider torus knots, $p$ and $q$ are relatively prime. The following is due to Moser \cite{Moser1971}:

\begin{proposition} $S^3_{\pm1/n}(T_{p,q})$ is Seifert fibered with base orbifold $S^2(p,q,pqn\mp1)$ for $n>0$. \end{proposition}
\begin{proof}Let $M=S^3\smallsetminus\nu(T_{p,q})$ so that $M(\alpha)=S^3_{r/s}(T_{p,q})$ for a given slope $\alpha=r\mu+s\lambda$, where $\mu$ is the knot meridian and $\lambda$ is the preferred longitude. As the complement of a regular fibre of a Seifert fibration of $S^3$, $M$ is Seifert fibered with base orbifold $D^2(p,q)$. Let $\fibre$ denote a regular fibre in $\partial M$; it is well known that $\fibre=pq\mu+\lambda$ (see Moser \cite{Moser1971}, for example). 

Now $M(\alpha)$ is Seifert fibered with base orbifold $S^2(p,q,\Delta(\alpha,\fibre))$ whenever $\alpha\ne\fibre$, according Heil \cite{Heil1974} (see also Moser \cite{Moser1971}). Here, $\Delta(\alpha,\fibre)$ measures the distance (i.e. the minimal geometric intersection number) between slopes the $\alpha$ and $\fibre$ in $\partial M$. 

In the present setting, $\alpha=\pm\mu+n\lambda$ for $n>0$ so that $M(\alpha)=S^3_{\pm1/n}(T_{p,q})$. Therefore, \[
\Delta(\alpha,\fibre)=|(\pm\mu+n\lambda)\cdot(pq\mu+\lambda)|=
\begin{cases}
pqn-1 & {\rm for\ positive\ surgeries} \\
pqn+1 & {\rm for\ negative\ surgeries}
\end{cases}\]
As a result, $M(\pm\mu+n\lambda)=S^3_{\pm1/n}(T_{p,q})$ is Seifert fibered with base orbifold $S^2(p,q,pqn\mp1)$ as claimed.
\end{proof}

\subsection*{Seifert involutions}

For a link $L\into S^3$, let $\Br(S^3,L)$ denote the two-fold branched cover of $S^3$, branched over $L$. The following is due to Seifert \cite{Seifert1933}:

\begin{proposition}$S^3_{\pm1/n}(T_{2,q})\cong\Br(S^3,T_{q,2qn\mp1})$ for $n>0$ and odd $q>1$. \end{proposition}
\begin{proof} The manifold $\Br(S^3,T_{q,2qn\mp1})$ is the Brieskorn sphere $\Sigma(2,q,2qn\mp1)$; this manifold is Seifert fibered with base orbifold $S^2(2,q,2qn\mp1)$ \cite[Lemma 1.1]{Milnor1975} (see also Seifert \cite[Zusatz zu Satz 17]{Seifert1933}). Since there is a  unique $\bZ$-homology sphere, for each $n>0$ and odd $q>1$, admitting a Seifert fibered structure with base orbifold $S^2(2,q,2qn\mp1)$ (see Scott \cite{Scott1983}, for example), the result follows.     \end{proof}

\subsection*{Montesinos involutions}

By a result of Schreier, the knot $T_{p,q}$ is strongly invertible \cite{Schreier1924}. As such, it is possible to realize the manifold $S^3_{r/s}(T_{p,q})$ as a two-fold branched cover via the Montesinos trick \cite{Montesinos1975}. We will adhere to the notation introduced in \cite[Section 3]{Watson2008} in constructing the relevant branch sets.

\begin{figure}[ht!]
\begin{center}
\raisebox{0pt}{\includegraphics[scale=0.4]{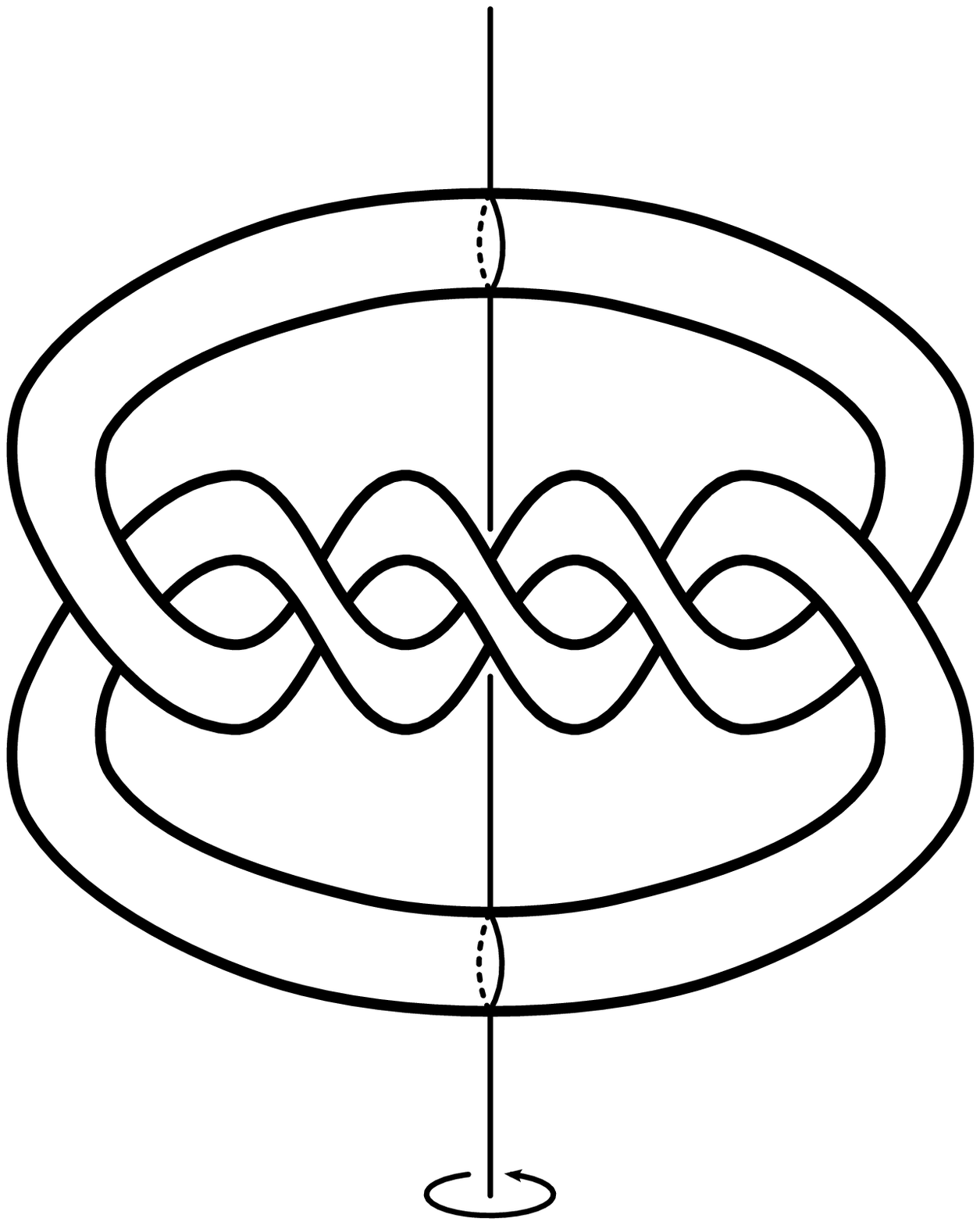}}\qquad
\raisebox{-40pt}{\includegraphics[scale=0.4]{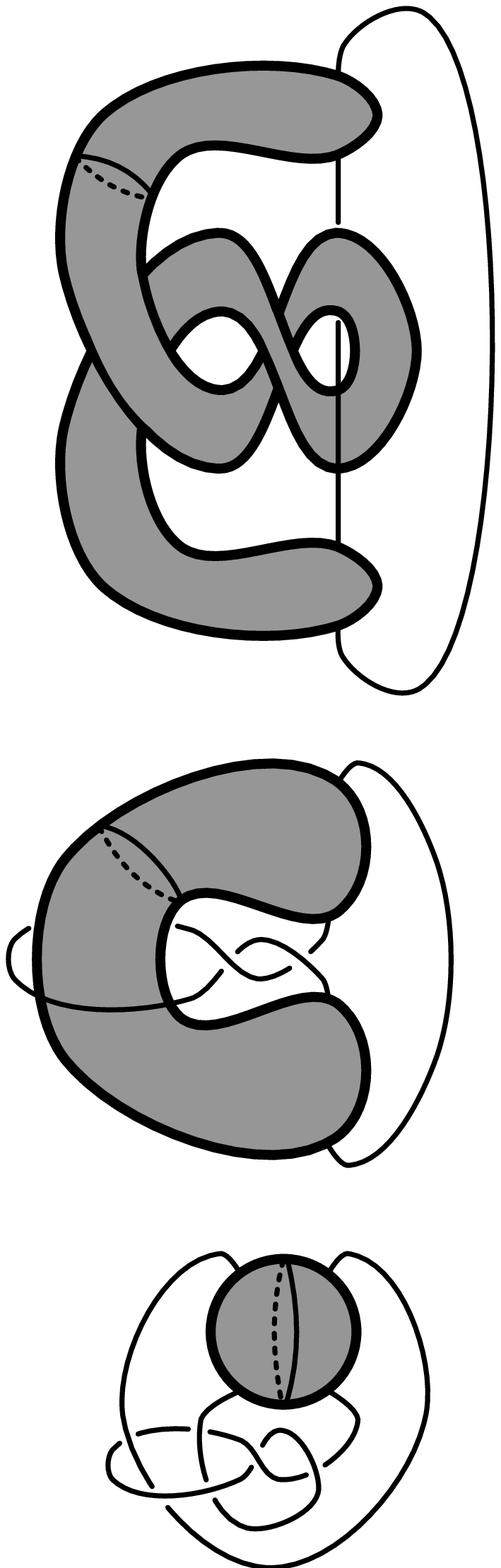}}\qquad
\labellist
	\pinlabel \rotatebox{-90}{$\cong$} at 163 429 
\endlabellist
\raisebox{0pt}{\includegraphics[scale=0.4]{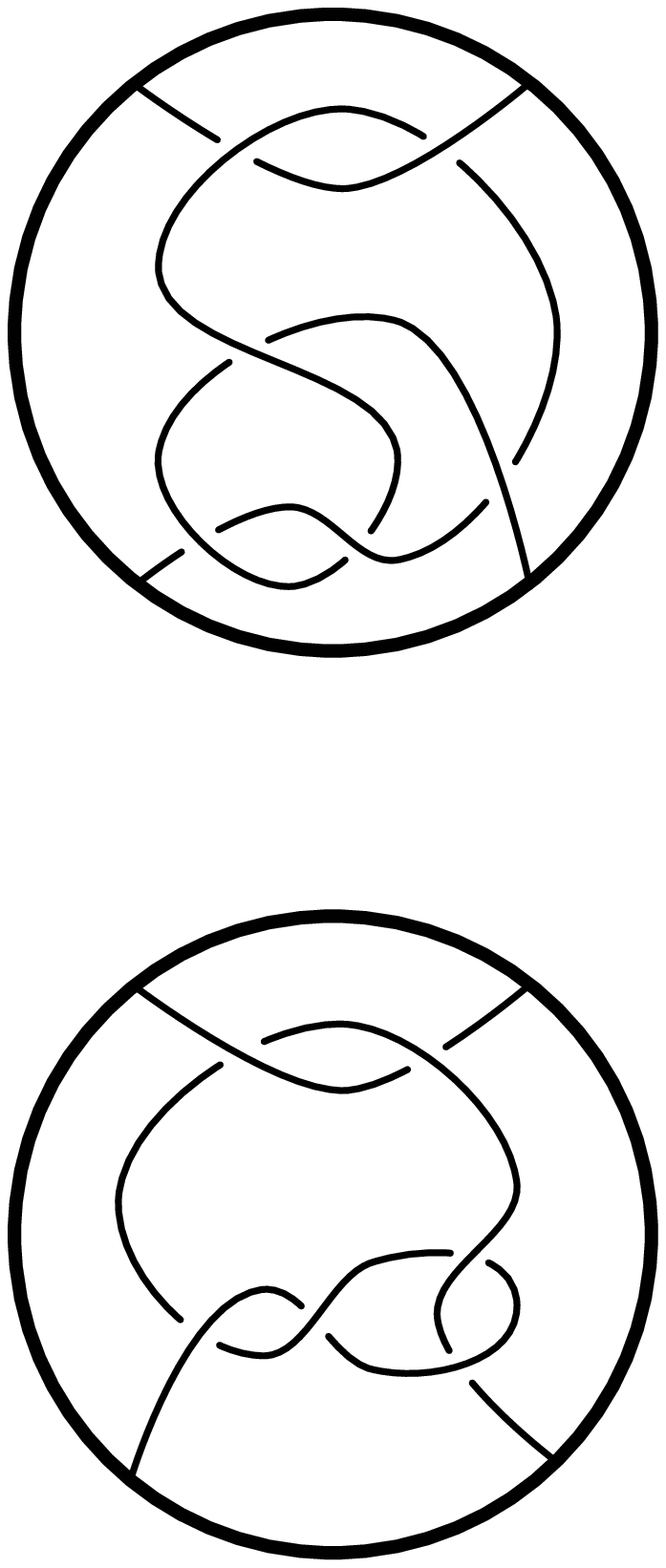}}
\end{center}
\caption{A strong inversion on the cinqfoil (left); isotopy of a fundamental domain (centre); and two representatives of the associated quotient tangle (right). Notice that the Seifert fibre structure on the knot complement is reflected in the sum of rational tangles of the associated quotient tangle.}
\label{fig:cinqfoil-inversion}\end{figure}

In the interest of being explicit, consider the cinqfoil $K=T_{2,5}$ (the knot $5_1$). A strong inversion on this knot is exhibited in Figure \ref{fig:cinqfoil-inversion}, together with an illustration of the process of obtaining a tangle with the property that $S^3\smallsetminus\nu(K)\cong\Br(B^3,\tau)$. That is, the complement of $K$ may be realized as the two-fold branched cover of a tangle $T=(B^3,\tau)$, where $\tau$ is the image of the fixed point set in the quotient.    

\begin{figure}[ht!]
\begin{center}
\raisebox{120pt}{\includegraphics[scale=0.4]{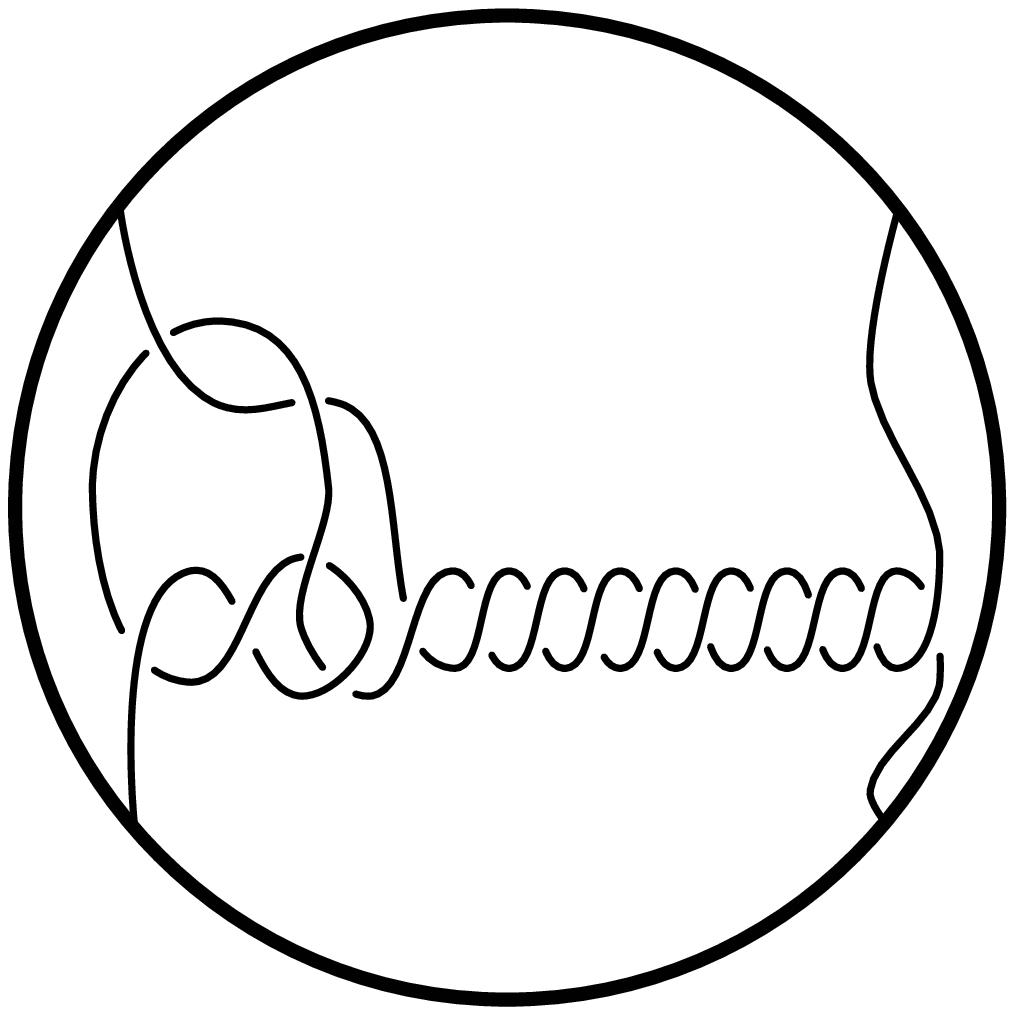}}\qquad\qquad
\labellist
	\small
	\pinlabel $\tau(\frac{1}{2})$ at 390 670
	\pinlabel $\tau(1)$ at 390 445
	\pinlabel $\tau(-1)$ at 395 209
	\pinlabel $\tau(-\frac{1}{2})$ at 395 -5
\endlabellist
\raisebox{0pt}{\includegraphics[scale=0.4]{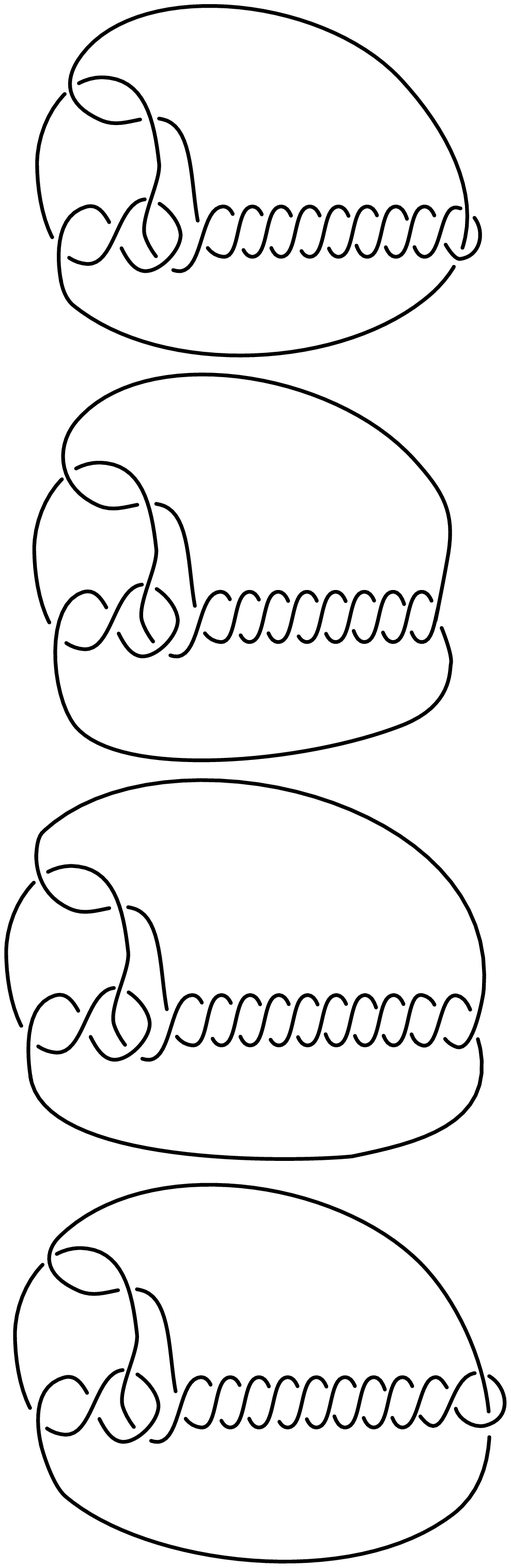}}
\end{center}
\caption{The preferred representative of the associated quotient tangle for the cinqfoil (left); and the branch sets $\tau(-\frac{1}{2})$, $\tau(-1)$, $\tau(1)$ and $\tau(\frac{1}{2})$ associated to $\{-\frac{1}{2},-1,1,\frac{1}{2}\}$-surgery on the cinqfoil (right).}
\label{fig:cinqfoil-framedtangle}\end{figure}

Note that such tangles are considered up to homeomorphism of the pair $(B^3,\tau)$ that need not fix the boundary in general. As a result, we may fix a preferred representative of such a tangle with the properties that 
\begin{itemize}
\item[(1)] the denominator closure of the tangle, denoted $\tau(\overzero)$, is unknotted and corresponds to a branch set for the trivial surgery, and \\
\item[(2)] the numerator closure, denoted  $\tau(0)$, gives a branch set for the zero surgery: $S^3_0(K)\cong\Br(S^3,\tau(0))$. 
\end{itemize}
This representative is illustrated in Figure \ref{fig:cinqfoil-framedtangle}. Note that it suffices to verify that $\det(\tau(0))=0$ to see that this is the appropriate framing on the tangle: this ensures that $\Br(S^3,\tau(0))$ has positive first Betti number. 

More generally, we have that $S^3_{r/s}(K)\cong\Br(S^3,\tau(\rs))$, where $\tau(\rs)$ is the link obtained by attaching an appropriate rational tangle. Note that the branch sets $\tau(n)$ associated to integer surgeries are obtained by adding $n$ half-twists (lifting to meridional Dehn twists in the cover). We illustrate $\tau(\pm1)$ and $\tau(\pm\frac{1}{2})$ in Figure \ref{fig:cinqfoil-framedtangle}, and refer the reader to \cite[Section 3]{Watson2008} for details.

We point to Montesinos' notes for a detailed discussion on Seifert fibered spaces as two-fold branched covers of $S^3$ in general \cite{Montesinos1976}.

\subsection*{Proof of the Theorem}

Continuing with $K=T_{2,5}$, by the observations above pertaining to the Seifert and Montesinos involutions we have that \[S^3_{\pm1/n}(K)\cong\Br(S^3,T_{5,10n\mp1})\cong\Br(S^3,\tau(\pm\textstyle\frac{1}{n}))\] for $n>0$. When $n=1$, using the program {\tt JavaKh} \cite{JavaKh} we calculate \[\rk\Khred(T_{5,10\mp1})=65\mp8 \ne 16\mp1=\rk\Khred(\tau(\pm1)).\] Similarly, when $n=2$ we calculate \[\rk\Khred(T_{5,20\mp1})=257\mp16 \ne 32\mp1=\rk\Khred(\tau(\pm\textstyle\frac{1}{2})).\] 

Each of these four pairs of examples illustrates a given manifold as a two-fold branched cover of $S^3$ in two different ways, with branch sets distinguished by the total rank of the reduced Khovanov homology. This proves the claim: $\rk\Khred$ is not an invariant of two-fold branched covers. 

\subsection*{Further remarks} We continue with the above notation for the tangle associated  to the cinqfoil.

\begin{proposition} $\rk\Khred(\tau(\pm\frac{1}{n}))\le 16n\mp1$ for $n>0$. \end{proposition}
\begin{proof}[Sketch of proof] Noting first that $\rk\Khred(\tau(\pm1))=16\mp1$, and calculating that $\rk\Khred(\tau(0))=16$, the result follows by induction on $n$: applying the long exact sequence for Khovanov homology we have that \[\textstyle\rk\Khred(\tau(\frac{1}{n})\le\rk\Khred(\tau(\frac{1}{n-1}))+\rk\Khred(\tau(0))=\rk\Khred(\tau(\frac{1}{n-1}))+16\] and \[\textstyle\rk\Khred(\tau(-\frac{1}{n})\le\rk\Khred(\tau(-\frac{1}{n-1}))+\rk\Khred(\tau(0))=\rk\Khred(\tau(-\frac{1}{n-1}))+16.\proved\]\end{proof}

While calculations of Khovanov homology for large torus knots are difficult to obtain, the existing calculations suggest that $\rk\Khred(T_{p,q})$ grows {\em at least} linearly in $q$. In particular, it seems reasonable to guess that surgery on the cinqfoil provides an infinite family of examples proving the Theorem. 

More generally, it would be interesting to understand the behaviour of the Khovanov homology for branch sets associated to $\frac{1}{n}$-framed surgery on the torus knots $T_{2,q}$ for $q\ge5$.

\bibliographystyle{plain}
\bibliography{/Users/liam/Documents/Mathematics/Bibliography/bibliography}

\begin{thebibliography}{10}

\bibitem{JavaKh}
Dror Bar-Natan and Jeremy Green.
\newblock {\tt JavaKh}.
\newblock Available at \url{http://www.katlas.org/wiki/KnotTheory}.

\bibitem{Bloom2009}
Jonathan Bloom.
\newblock {Odd Khovanov homology is mutation invariant}.
\newblock Preprint arXiv:0903.3746.

\bibitem{Heil1974}
Wolfgang Heil.
\newblock Elementary surgery on {S}eifert fiber spaces.
\newblock {\em Yokohama Math. J.}, 22:135--139, 1974.

\bibitem{Khovanov2000}
Mikhail Khovanov.
\newblock A categorification of the {J}ones polynomial.
\newblock {\em Duke Math. J.}, 101(3):359--426, 2000.

\bibitem{Khovanov2003}
Mikhail Khovanov.
\newblock Patterns in knot cohomology. {I}.
\newblock {\em Experiment. Math.}, 12(3):365--374, 2003.

\bibitem{Lipshitz2008}
Robert Lipshitz.
\newblock {Is Khovanov homology an invariant of the the branched double cover?}
\newblock From ``the {F}loer homology problem collection", available at
  \url{http://www.floerhomology.com}.

\bibitem{Milnor1975}
John Milnor.
\newblock On the {$3$}-dimensional {B}rieskorn manifolds {$M(p,q,r)$}.
\newblock In {\em Knots, groups, and 3-manifolds ({P}apers dedicated to the
  memory of {R}. {H}. {F}ox)}, pages 175--225. Ann. of Math. Studies, No. 84.
  Princeton Univ. Press, Princeton, N. J., 1975.

\bibitem{Montesinos1976}
Jos{\'e}~M. Montesinos.
\newblock {Rev\^etements ramifi\'es de n{\oe}ds, espaces fibr\'e de Seifert et
  scindements de Heegaard}.
\newblock Lecture notes, Orsay 1976.

\bibitem{Montesinos1975}
Jos{\'e}~M. Montesinos.
\newblock Surgery on links and double branched covers of {$S\sp{3}$}.
\newblock In {\em Knots, groups, and $3$-manifolds (Papers dedicated to the
  memory of R. H. Fox)}, pages 227--259. Ann. of Math. Studies, No. 84.
  Princeton Univ. Press, Princeton, N.J., 1975.

\bibitem{Moser1971}
Louise Moser.
\newblock Elementary surgery along a torus knot.
\newblock {\em Pacific J. Math.}, 38:737--745, 1971.

\bibitem{Ozsvath2008}
Peter Ozsv{\'a}th.
\newblock {Private communication}.

\bibitem{OSz2005-branch}
Peter Ozsv{\'a}th and Zolt{\'a}n Szab{\'o}.
\newblock On the {H}eegaard {F}loer homology of branched double-covers.
\newblock {\em Adv. Math.}, 194(1):1--33, 2005.

\bibitem{Schreier1924}
Otto Schreier.
\newblock {\"Uber die Gruppen $A^aB^b=1$}.
\newblock {\em Abh. Math. Sem. Univ. Hamburg}, 3:167--169, 1924.

\bibitem{Scott1983}
Peter Scott.
\newblock The geometries of {$3$}-manifolds.
\newblock {\em Bull. London Math. Soc.}, 15(5):401--487, 1983.

\bibitem{Seifert1933}
H.~Seifert.
\newblock Topologie {D}reidimensionaler {G}efaserter {R}\"aume.
\newblock {\em Acta Math.}, 60(1):147--238, 1933.

\bibitem{Watson2008}
Liam Watson.
\newblock Surgery obstructions from {K}hovanov homology.
\newblock Preprint arXiv:0807.1341.

\bibitem{Wehrli2009}
Stephan Wehrli.
\newblock {Mutation invariance of Khovanov homology over $\mathbb{F}_2$}.
\newblock Preprint arXiv:0904.3401.

\bibitem{Wehrli2007}
Stephan Wehrli.
\newblock {Mutation invariance of Khovanov homology over $\mathbb{Z}_2$}.
\newblock Lecture notes, Kyoto 2007.

\end{thebibliography}

\end{document}